

\input epsf.tex


\parindent=0.7truecm
\parskip=3pt plus .5pt minus .5pt
\hsize=16truecm \vsize=22.7truecm

\let\Fi=\varphi \let\eps=\varepsilon
\font\bigbf=cmbx10 scaled \magstephalf     
\font\Bigbf=cmbx12 scaled \magstep1        
\font\letranota=cmr9
\font\EU=eufm10


\font\TenEns=msbm10
\font\SevenEns=msbm7
\font\FiveEns=msbm5
\newfam\Ensfam
\def\Ens{\fam\Ensfam\TenEns}
\textfont\Ensfam=\TenEns
\scriptfont\Ensfam=\SevenEns
\scriptscriptfont\Ensfam=\FiveEns

\def\R{{\Ens R}}

\def\N{{\Ens N}}

\font\TenCM=cmr10
\font\SevenCM=cmr7
\font\FiveCM=cmr5
\newfam\CMfam
\def\CM{\fam\CMfam\TenCM}
\textfont\CMfam=\TenCM
\scriptfont\CMfam=\SevenCM
\scriptscriptfont\CMfam=\FiveCM

\font\TenVect=cmbxti10
\font\SevenVect=cmbxti7
\newfam\Vectfam
\def\Vect{\fam\Vectfam\TenVect}
\textfont\Vectfam=\TenVect
\scriptfont\Vectfam=\SevenVect
\scriptscriptfont\Vectfam=\SevenVect
  \def\letravetor#1{{\Vect #1\/}}

\font\TenEU=eufm10
\font\SevenEU=eufm7
\newfam\EUfam
\def\EU{\fam\EUfam\TenEU}
\textfont\EUfam=\TenEU
\scriptfont\EUfam=\SevenEU
\scriptscriptfont\EUfam=\SevenEU
   \def\letraEU#1{{\EU #1\/}}


\def\unim{\,{\CM i{}\/}}
\def\ee{\,{\CM e\/}}
\def\sign{\mathop{\rm sign}\nolimits}
\def\mn{{\hbox{\fiverm min}}}
\def\app{{\hbox{\fiverm app}}}

\let\Fi=\varphi
\let\eps=\varepsilon
\def\mod#1{\vert #1\vert}
\def\Mod#1{\left|#1\right|}
\def\norma#1#2{\Vert #1\Vert_{#2}}
\setbox111=\hbox to 2truemm{\hrulefill}
\setbox222=\hbox to 2truemm{\vrule height 2truemm width .4truept\hfil\vrule height 2truemm width .4truept}
\def\cqd{\vbox{\offinterlineskip\copy111\copy222\copy111}}
\def\build#1_#2^#3{\mathrel{\mathop{\kern 0pt#1}\limits_{#2}^{#3}}}


\newcount\numerosection
\newcount\eqnumer
\newcount\lemnumer
\newcount\fignumer
\numerosection=0
\eqnumer=0
\lemnumer=0
\fignumer=0
\def\numsection{\global\advance\numerosection by1
\global\eqnumer=0
\global\lemnumer=0
\global\fignumer=0
\the\numerosection}
\def\strutdepth{\dp\strutbox}
\def\marginalsigne#1{\strut
    \vadjust{\kern-\strutdepth\specialsigne{#1}}}
\def\specialsigne#1{\vtop to \strutdepth{
\baselineskip\strutdepth\vss\llap{#1 }\null}}
\font\margefont=cmr10 at 6pt
\newif\ifshowingMacros
\showingMacrosfalse

\def\cite#1{\csname#1\endcsname}
\def\label#1{\gdef\currentlabel{#1}}
\def\Lefteqlabel#1{\global\advance\eqnumer by 1
\label{#1}
\ifx\currentlabel\relax
\else
\expandafter\xdef
\csname\currentlabel\endcsname{(\the\numerosection.\the\eqnumer)}
\fi
\global\let\currentlabel\relax
\ifshowingMacros
\leqno\llap{%
\margefont #1\hphantom{M}}(\the\numerosection.\the\eqnumer)%
\else
\leqno(\the\numerosection.\the\eqnumer)%
\fi
}

\def\Righteqlabel#1{\global\advance\eqnumer by 1
\label{#1}
\ifx\currentlabel\relax
\else
\expandafter\xdef
\csname\currentlabel\endcsname{(\the\numerosection.\the\eqnumer)}%
\fi
\global\let\currentlabel\relax
\ifshowingMacros
\eqno(\the\numerosection.\the\eqnumer)%
\rlap{\margefont\hphantom{M}#1}
\else
\eqno(\the\numerosection.\the\eqnumer)%
\fi
}
\let\reqlabel=\Righteqlabel
\def\numlabel#1{\global\advance\eqnumer by1
\label{#1}
\ifx\currentlabel\relax
\else
\expandafter\xdef
\csname\currentlabel\endcsname{(\the\numerosection.\the\eqnumer)}
\fi
\global\let\currentlabel\relax
\ifshowingMacros
   (\the\numerosection.\the\eqnumer)\rlap{\margefont\hphantom{M}#1}
\else
(\the\numerosection.\the\eqnumer)%
\fi
}
\def\numero{\global\advance\eqnumer by1
\number\numerosection.\number\eqnumer}
\def\lemlabel#1{\global\advance\lemnumer by 1
\ifshowingMacros%
   \marginalsigne{\margefont #1}%
\else
    \relax
\fi
\label{#1}%
\ifx\currentlabel\relax%
\else
\expandafter\xdef%
\csname\currentlabel\endcsname{\the\numerosection.\the\lemnumer}%
\fi
\global\let\currentlabel\relax%
\the\numerosection.\the\lemnumer}
\def\numlem{\global\advance\lemnumer by1
\the\numerosection.\number\lemnumer}


\let\itemBibli=\item
\def\bibl#1#2\endbibl{\par{\itemBibli{#1} #2\par}}
\def\ref.#1.{{\csname#1\endcsname}}
\newcount\bib   \bib=0
\def\bibmac#1{\advance\bib by 1
\expandafter
\xdef\csname #1\endcsname{\the\bib}}
\def\BibMac#1{\advance\bib by1
\bibl{\letranota[\the\bib]}%
{\csname#1\endcsname}\endbibl}
\def\MakeBibliography#1{
\noindent{\bf #1}
\bigskip
\bib=0
\let\bibmac=\BibMac
\BiblioFil
\BiblioOrd
}

\newif\ifbouquin

\bouquinfalse
\font\tenssi=cmssi10 at 10 true pt
 at 10 true pt
\def\bibliostyle#1#2#3{{\rm #1}\ifbouquin{\tenssi #2\/}\global\bouquinfalse\else{\it #2\/}\fi{\rm #3}}


\def\BiblioFil{
  \def\Bose{\bibliostyle{\letranota S.N.~Bose, Z.~Phys., { 26}, 178, 1924.}{}{}}
  \def\EinsteinA{\bibliostyle{\letranota A.~Einstein, Sitzungsber.~K.~Preuss.~Akad.~Wiss.~Phys.~Math., { 261}, (1924).}{}{}}
  \def\EinsteinB{\bibliostyle{\letranota A.~Einstein, Sitzungsber.~K.~Preuss.~Akad.~Wiss.~Phys.~Math., { 3}, (1925).}{}{}}
  \def\GSS{\bibliostyle{\letranota A.~Griffin, D.W.~Snoke and Stringari (ed), Bose-Einstein Condensation,
       Cambridge University Press, (1995).}{}{}}
  \def\AEMW{\bibliostyle{\letranota M.H.J.~Anderson, J.R.~Ensher, M.R.~Matthews and
       C.E.~Wiemna, Science, { 269}, (1995), 198.}{}{}}
  \def\DMADDKK{\bibliostyle{\letranota K.B.~Davis at al., Phys.~Rev.~Lett., { 75}, (1995), 3969.}{}{}}
  \def\Gross{\bibliostyle{\letranota E.P.~Gross, J.~Math.~Phys., { 4}, (1963), 195.}{}{}}
  \def\Pitaev{\bibliostyle{\letranota L.K.~Pitaevskii, Sov.~Phys.~JETP, { 13}, (1961), 451.}{}{}}
  \def\Carretero{\bibliostyle{\letranota R.~Carretero-Gonz\'alez, D.J.~Frantzeskakis and P.G.~Kevrekidis,
       Nonlinearity, { 21}, (2008), R139.}{}{}}
  \def\LiebEtAl{\bibliostyle{\letranota E.H.~Lieb, R.~Seiringer and J.~Yngvason, Physical Review A, { 61}, (2000), 043602.}{}{}}
  \def\MunozDelgado{\bibliostyle{\letranota A.M.~Mateo and V.~Delgado, Phys.~Rev.~E {88}, (2013), 042916.}{}{}}
  \def\KavWeis{\bibliostyle{\letranota O.~Kavian and F.~Weissler, Mich.~J.~Math., 41, (1994), 151.}{}{}}
  \def\LiebLoss{\bibliostyle{\letranota E.H.~Lieb and M.~Loss, Analysis 2nd Edition, Graduate Studies in
       Mathematics, Vol.~14, AMS, (2001).}{}{}}
  \def\Mossino{\bibliostyle{\letranota J.~Mossino, Inegalit\'es Isoperimetriques et Aplications en Physique, Hermann, (1984).}{}{}}
  \def\Carles{\bibliostyle{\letranota R.~Carles, Ann.~Henri Poincar\'e, { 3}, (2003), 757.}{}{}}
  \def\Oh{\bibliostyle{\letranota Y-G.~Oh, J.~Diff.~Eq., {81}, (1989), 255.}{}{}}
  \def\CazeLions{\bibliostyle{\letranota T.~Cazenave T and P-L.~Lions, Comm.~Math.~Phys., { 85}, (1982), 153.}{}{}}
  \def\CipGondTrall{\bibliostyle{\letranota R.~Cipolatti, J.~L\'opez Gondar~and C.~Trallero-Giner, Physica D, {241}, (2012), 755.}{}{}}
  \def\TrallCipLiew{\bibliostyle{\letranota C.~Trallero-Giner, R.~Cipolatti and T.C.H.~Liew , Eur.~Phys.~J.~D, {67}, (2013), 143.}{}{}}
}

\def\BiblioOrd{
\bibmac{MunozDelgado}    
\bibmac{EinsteinA}            
\bibmac{EinsteinB}            
\bibmac{Bose}                  
\bibmac{GSS}                   
\bibmac{AEMW}                
\bibmac{DMADDKK}          
\bibmac{Gross}                 
\bibmac{Pitaev}                
\bibmac{Carretero}           
\bibmac{LiebEtAl}             
\bibmac{KavWeis}            
\bibmac{LiebLoss}            
\bibmac{Mossino}             
\bibmac{Carles}               
\bibmac{Oh}                    
\bibmac{CazeLions}         
\bibmac{CipGondTrall}     
\bibmac{TrallCipLiew}      
}

\BiblioOrd



\centerline{\Bigbf Mathematical analysis of a Mu\~noz-Delgado model }
\centerline{\Bigbf for cigar-shaped Bose-Einstein condensates}

\bigskip
\centerline{R.~Cipolatti$^{1}$\footnote{*}{\letranota Corresponding author; E-mail address: cipolatti@im.ufrj.br.},
      Y.M.~Lira$^{2}$ and G.~Saisse$^{3}$}
\bigskip

{\baselineskip=14truept
\letranota
\centerline{$^{1}\,$Instituto de Matem\'atica, Universidade Federal do Rio de Janeiro}
\centerline{Ilha do Fund\~ao, Rio de Janeiro, RJ, Brasil, CEP: 21941-909}
\centerline{$^{2}\,$Instituto de Matem\'atica e Estat{\'\i}stica,  Universidade do Estado  do Rio de Janeiro}
\centerline{Maracan\~a, Rio de Janeiro, RJ, Brasil, CEP 20550-013}
\centerline{$^{3}\,$Instituto de F{\'\i}sica, Universidade Federal do Rio de Janeiro}
\centerline{Ilha do Fund\~ao, Rio de Janeiro, RJ, Brasil, CEP 21941-909}

}


\bigskip
\hrule
\bigskip

{\letranota\baselineskip=12truept
\noindent{\bf Abstract} In this paper we present mathematical analysis of one-dimensional effective models
proposed in [\cite{MunozDelgado}] concerning Bose-Einstein condensates in the presence of harmonic confinement.
Among the demonstrated properties,  we can mention: existence, uniqueness, orbital stability, symmetry
and gaussian asymptotic decay of ground-state solutions in the repulsive case.
We also report formul\ae\ for the minimal energy $E_{\mn}$  and the associate chemical potential $\mu$
as functions of a parameter $\lambda$, which is related to $N$ (the number of atoms) and/or $a$
(the s-wave scattering length). By considering Taylor's development of the non-quadratic therm of the
energy and using appropriate gaussian functions as approximations for the ground state,
we present some numerical experiments to illustrate our results.

\medskip
\noindent{\letranota Key words: Bose-Einstein condensates, stability of ground states, analytical
   approximate formulae, repulsive or attractive interatomic interactions.}
}

\bigskip
\hrule
\bigskip  
\noindent{\bigbf\numsection.\ Introduction}\par
\smallskip  

\noindent
The phenomenon of Bose-Einstein condensates, predicted by Einstein in 1924 [\cite{EinsteinA},\cite{EinsteinB}] based on the statistics
proposed by Bose [\cite{Bose}] for boson gases, are connected to superfluidity in liquid helium and superconductivity in metals [\cite{GSS}].
They were experimentally realized seventy years after their prediction for different species of gases confined in  magnetic and/or optical
traps under extremely low temperatures [\cite{AEMW},\cite{DMADDKK}]. For these reasons, this phenomenon has attracted
the attention of many scientists, specially during recent years.

From a theoretical point of view, the dynamics of a dilute ultracold atomic gas composed by $N$ interacting bosons of mass $m$
confined by an external potential $V_{\hbox{\sevenrm ext}}(\letravetor x)$ can be described by the Heisenberg evolution equation
$$\unim\hbar{\partial\Psi\over\partial t}(t,\letravetor x) =  \left[-{\hbar^2\over 2m}\Delta+V_{\hbox{\sevenrm ext}}(\letravetor x)
    +\int_{\R^3}\overline{\Psi(t,\letravetor y)}V(\letravetor x-\letravetor y)\Psi(t,\letravetor y)\,d\letravetor y\right]\Psi(t,\letravetor x),
    \reqlabel{Heisen}$$
where $\hbar$ is Dirac's constant, $\Psi$ is the macroscopic wave function of the condensate (the expetaction value of the
field operator),  $V(\letravetor x-\letravetor y)$ is the interatomic potential (the  two-body interaction field) and $\letravetor{x}=(x,y,z)$.

Under near absolute zero temperatures, where thermal excitations can be neglected, the macroscopic wave function $\Psi(t,\letravetor x)$
can be accurately described by the three dimensional (3D) Gross-Pitaevskii equation (GPE) [\cite{Gross},\cite{Pitaev}], namely
$$\unim\hbar{\partial\Psi\over\partial t}(t,\letravetor x) =  \left[-{\hbar^2\over 2m}\Delta+V_{\hbox{\sevenrm ext}}(\letravetor x)
    +gN\mod{\Psi(t,\letravetor x)}^2\right]\Psi(t,\letravetor x),\reqlabel{GPE}$$
where $g=4\pi\hbar^2a/m$ is the interection strength determined by the $s$-wave scattering length $a$  ($a<0$ when interatomic
forces are attractive and $a>0$ for repulsive interatomic forces), for which the BECs are characterized by its ground-state solutions.

It is usual to apply external potentials $V_{\hbox{\sevenrm ext}}(\letravetor x)$ to trap  the condensate.
In the case of magnetic ones, the external potential has the harmonic form
$V_{\hbox{\sevenrm ext}}(\letravetor x):={1\over 2}m\letravetor{w}\!\cdot\!\letravetor{x}$,
where $\letravetor{w}=(\omega_1,\omega_2,\omega_3)$, with which one can manipulate the condensate
using different frequencies along the three directions.  The flexibility over the choice of the confining
frequencies may be used to control its shape. For example, the strongly anisotropic case with
$\omega_\perp:=\omega_2=\omega_3\gg\omega_1$ is particularly interesting as it is related
to effective quasi-one-dimensional (1D) BEC.

Under certain choices of the physical parameters (see [\cite{Carretero}]), the transverse confinement
of the condensate is so tight  that the dynamics of such a cigar-shaped BEC can be considered to be effectively 1D.
However, this dimensional reduction should only be considered as a 1D limit of a 3D mean-field
theory, instead of a genuine 1D model [\cite{LiebEtAl}].  Moreover, since the numerical resolution of \cite{GPE}
is a heavy task that requires considerable computational effort, several authors have looked for
more effective one-dimensional models to analyze the behavior of cigar-shaped BECs.

This is the case for a one-dimensional model derived by  A.~Mu\~noz Mateo and V.~Delgado,
(see [\cite{MunozDelgado}] for more details). They assume that the axial degrees of freedom evolve so slowly in
time in comparison with the transverse degrees of freedom in such a way that, at each instant $t$, correlations
between axial and radial motions are negligible and the wave function can be factorized such that
$\Psi(t,\letravetor{x}) = \Phi(y,z,n_1)\phi(t,x)$, where $n_1$ is the axial linear density,
$$ n_1(t,x) := N\int_{\R^2}\mod{\Phi(t,x,y,z)}^2\,dydz=N\mod{\phi(t,x)}^2$$
with both $\Phi$ and $\phi$ normalized to unity.
This procedure allows us to decompose the equation \cite{GPE} in such a way that
$$  \unim\hbar{\partial\phi\over \partial t}=-{\hbar^2\over 2m}{\partial^2\phi\over\partial x^2}
     +{1\over 2m}\omega_1^2x^2 +\mu_\perp(n_1)\phi,$$
where $\mu_\perp(n_1)$ is the local chemical potential satisfying the stationary transverse  GPE
$$ \left[-{\hbar^2\over 2m}\left({\partial^2\over\partial y^2}+{\partial^2\over\partial z^2}\right) + {1\over 2m}\left(\omega_2^2y^2
  +\omega_3^2z^2\right)+gn_1\mod{\Phi}^2\right]\Phi=\mu_\perp(n_1)\Phi.$$

After a cumbersome procedure using decomposition of generalized Laguerre polynomials and physical considerations,
they have derived that
$$ \mu_\perp(n_1)=\hbar\omega_\perp\alpha{1+3\gamma a n_1 \over \sqrt{1+2\gamma a n_1}},$$
which leads us to the following equation for the axial dynamics in the $x$ direction,
$$ \unim\hbar{\partial\phi\over\partial t} = -{\hbar^2\over 2m}{\partial^2\phi\over \partial x^2}+{1\over 2m}\omega_1^2x^2\phi
    +\hbar\omega_\perp\alpha{1+3\gamma aN\mod{\phi}^2   \over  \sqrt{1+2\gamma a N\mod{\phi}^2}} \phi. \reqlabel{MDE}$$

In this paper we are interested in the mathematical analysis of ground-state solutions of equation \cite{MDE}, {\it i.e.},
the standing wave solutions of minimal energy: $\phi(t,x)=\phi_0(x)\ee^{-\unim\mu_0 t/\hbar}$, which leads us to consider the equation
$$ -{\hbar^2\over 2m}{d^2\phi_0\over dx^2}+{1\over 2m}\omega_3^2x^2\phi_0
    +\hbar\omega_\perp\alpha{1+3\gamma aN\mod{\phi_0}^2   \over  \sqrt{1+2\gamma a N\mod{\phi_0}^2}} \phi_0
    = \mu_0\phi_0. \reqlabel{StatMDE}$$

To proceed in this goal, we consider the dimensionless formulation of the previous equation, which takes the form
$$ -{d^2\Fi\over ds^2} + s^2\Fi +C_\omega{1+3\lambda\mod{\Fi}^2 \over \sqrt{1+2\lambda\mod{\Fi}^2}} \Fi
    = \mu\Fi, \reqlabel{OurMDE}$$
where,  for $ l_0:=\sqrt{\hbar/m\omega_1}$, we set $s:=x/l_0$, $\Fi(s):=\phi_0(l_0s)$ and
$$ \mu:= {\mu_0\over \hbar\omega_1},\quad C_\omega:= {2\omega_\perp\alpha\over \omega_1},\quad
    \lambda:={\gamma a N\over l_0}.$$

Note that the solutions of \cite{OurMDE} can be viewed as standing waves of the following dimensionless time dependent GPE
$$ \unim{\partial u\over \partial\tau}=-{\partial^2 u\over\partial s^2} + s^2u
   +C_\omega{1+3\lambda\mod{u}^2 \over \sqrt{1+2\lambda\mod{u}^2}}u. \reqlabel{TimeDepOurMDE}$$

We organize the paper as follows: in Section~2 we prove that every complex solution of equation \cite{OurMDE} has the form
$\psi(s)=\ee^{\!\unim\theta}\Fi(s)$ were $\theta\in\R$ and $\Fi$ is a real solution that decays at infinite as gaussian.
We prove the existence and uniqueness of ground states for every $\lambda\ge0$ and that these solutions
are orbitally stable. In Section~3 we deduce a general formula relating the minimal energy and the associated chemical
potential as functions of the parameter $\lambda$. In Section~4 we present some numerical results and, finally, in Section~5,
we present an extension for the present context of Thomas-Fermi's approximation method.

\medskip 
\noindent{\bf\numsection.\ Existence and stability of ground states in the repulsive case ($\lambda>0$)}\par
\smallskip 

Although the solutions of \cite{TimeDepOurMDE} are in general complex valued functions and because we are in dimension one,
Sobolev's embedding theorems alows  us to restrict our analysis of the existence of ground states for real valued ones,
as we can see by the following lemma, where $H^1(\R)$ denotes the usual Sobolev space.

\smallskip
\noindent{\bf Lemma\ \lemlabel{Lemma1-1}:} {\sl If $\psi\in H^1(\R)$ is a complex solution of \cite{OurMDE} (which means
$\psi=u+\unim v$ whith $u,v\in H^1(\R)$), then there exists a real function
$U(s)$ solution of $\cite{OurMDE}$ and a real number $\theta$ such that $\psi(s)=\ee{}\!^{\unim\theta}U(s)$.}

\noindent
{\bf Proof:}  To simplify the notation, we define $f_\lambda:[0,+\infty)\rightarrow\R$ as
$$ f_\lambda(\rho):=  {1+3\lambda\rho^2  \over \sqrt{1+2\lambda\rho^2} }
      = {3\over 2}\sqrt{1+2\lambda\rho^2} - {1\over 2}{1\over\sqrt{1+2\lambda\rho^2}}. \reqlabel{DefFuncf}$$
Then the equation \cite{OurMDE} reads as
$$ -{d^2\hfil\over ds^2}\Fi(s) + s^2\Fi(s) +C_\omega f_\lambda\bigl(\mod{\Fi(s)}\bigr) \Fi(s)  = \mu\Fi(s),\quad s\in\R. \reqlabel{OurMDEbis}$$
If $\psi=u+\unim v$ is a complex solution of \cite{OurMDEbis},  it follows that
$$\left\{\eqalign{
-{d^2 u \over ds^2}+s^2 u+ C_\omega f_\lambda\bigl(\mod{\psi}\bigr)u & =\mu u,\cr
-{d^2 v \over ds^2}+s^2 v+ C_\omega f_\lambda\bigl(\mod{\psi}\bigr)v & =\mu v.\cr}
\right.\reqlabel{SistRealSol}$$
Now, multiplying the first equation in the above system by $v$ and the second by $u$, we get
$$u{d^2v\over ds^2}-v{d^2u\over ds^2}=0\quad\Rightarrow\quad u{dv\over ds }-v{du\over ds}=C,\quad \forall s\in\R,$$
for some real constant $C$. As it is well known that $H^1(\R)\subset C_0(\R)$, where $C_0(\R)$ is the space of continuous function
which tend to zero as $s\rightarrow\pm\infty$, we conclude that $C=0$ and hence
$$u{dv\over ds}-v{du\over ds }=0\quad\Rightarrow\quad {d\hfil\over ds}\left({v\over u}\right)=0.$$
Therefore,  $u=\beta v$, $\beta\in\R$ and each one of the equations of \cite{SistRealSol} reduces to
$$-{d^2 u \over ds^2}+s^2 u+C_\omega f_\lambda\bigl((1+\beta^2)^{1/2}\mod{u}\bigr)u =\mu u.$$
Now, considering $U(s):=(1+\beta^2)^{1/2}u(s)$, it follows that $U(s)$ is a real valued solution of \cite{OurMDEbis} and
$$\psi=u+\unim v=\left({1\over\sqrt{1+\beta^2}}+{\unim\beta\over\sqrt{1+\beta^2}}\right)U=\ee{}\!^{\unim\theta}U,$$
where $\theta:=\arctan\beta$.\quad\cqd

In the sequel we denote $\norma{\psi}{p}$ the usual norm of $\psi\in L^p(\R)$.

\smallskip
\noindent
{\bf Lemma\ \lemlabel{Lemma2}:} {\sl If $\psi\in H^1(\R)$ is a real solution of \cite{OurMDE}, then $\mu>1+C_\omega$.}

\noindent
{\bf Proof:}  Indeed, let $\varphi_1:\R\rightarrow\R$ be the function
$$\varphi_1(s):={1\over\root4\of\pi}\exp(-s^2/2).\reqlabel{Psi1}$$
It is easy to see that $\norma{\varphi_1}{2}=1$ and that $-\varphi_1''(s)+s^2\varphi_1(s)=\varphi_1(s)$. This means that $\varphi_1$ is
an eigenfunction of the operator $L:=-{d^2\hfil\over ds^2}+s^2$ corresponding to the eigenvalue $\lambda_1=1$. In fact, $L$ has an
infinite sequence of eigenvalues $\lambda_1<\lambda_2<\cdots$, where $\lambda_n=(2n-1)$,  $n\in\N$, and the
Hermite functions are the corresponding eigenfunctions. It is also known that $\lambda_1$ has the following variational characterization,
$$\lambda_1=\min\Bigl\{\int_\R\left(\Mod{\psi'(s)}^2+s^2\mod{\psi(s)}^2\right)ds \,;\,
       \norma{\psi}{2}=1\Bigr\}.\reqlabel{DefLambda1}$$

So, if we multiply \cite{OurMDE} by $\Fi$ and integrate on $\R$, we get
$$\mu\norma{\Fi}{2}^2 = \int_\R\left(\Mod{\Fi'(s)}^2
    +s^2\mod{\Fi(s)}^2+C_\omega f_\lambda\bigl(\mod{\Fi(s)}\bigr)\mod{\Fi(s)}^2)\right)ds.$$
Since $f_\lambda(\rho)\ge1$ for all $\rho\ge 0$, we get from \cite{DefLambda1}
$$ \mu\norma{\Fi}{2}^2\ge \bigl(\lambda_1+C_\omega)\bigr)\norma{\Fi}{2}^2\quad\Rightarrow
   \quad \mu\ge 1+C_\omega.\quad\cqd$$

An important property of solutions of \cite{OurMDE} is their gaussian asymptotic decay at infinity, as asserted in the following result.

\smallskip
\noindent {\bf Theorem\ \lemlabel{ExpDecay}:} {\sl Let  $\Fi\in H^1(\R)$ be a solution of \cite{OurMDE}.
Then $\Fi\in C^2(\R)$ and there exists $k_0\in (0,1)$ and $C(k_0)>0$ such that
$$\mod{\Fi(s)}\le C(k_0)\exp[-k_0s^2/2]\quad \forall s \in\R.\reqlabel{DecaiExp}.$$}\relax
{\bf Proof:}
We proceed as in [\cite{KavWeis}].  Since $\Fi\in H^1(\R)$, it follows from Sobolev embedding that $\Fi(s)$ is a continuous
function satisfying
$$\lim_{\mod{s}\to+\infty}\varphi(s)=0.\reqlabel{Tend2Zero}$$
If we denote $g(s):=s^2-\mu+C_\omega f_\lambda \bigl(\mod{\Fi(s)}\bigr)$, it follows from \cite{OurMDE},
$$ \bigl<\Fi'':\phi\bigr> = \bigl<g\Fi:\phi\bigr>, \quad \forall\phi\in D(\R). \reqlabel{OurMDE4Kato}$$
Let $h(s):=\mod{\Fi(s)}$. By Kato's inequality we have $h''\ge\sign(\Fi)\Fi''$ in the sense of distributions.
So, by multiplying \cite{OurMDE4Kato} by $\sign(\Fi)$, we get
$$ \bigl<h'':\phi\bigr> \ge  \bigl<\sign(\Fi)\Fi'':\phi\bigr> =  \bigl<g\mod{\Fi}:\phi\bigr>,
    \quad \forall\phi\in D(\R),\quad\phi\ge0,$$
which means that
$$ -h''+gh\le 0 \hbox{\rm\ \  in the sense of distributions on\ } \R.\reqlabel{PrinMax1}$$

On the other hand, if we set $\Fi_\kappa(s):=\!\ee^{-\kappa s^2/2}$,   a simple calculation gives
$$-\Fi_\kappa''(s)+g(s)\Fi_\kappa(s) = \bigl[(1-\kappa^2)s^2+1-\mu+C_\omega f_\lambda\bigl(\mod{\Fi(s)}\bigr)\bigr]\Fi_\kappa(s).$$
Since $f_\lambda$ is a positive increasing function on $[0,+\infty)$ and $f_\lambda(0)=1$, it follows that
$$-\Fi_\kappa''(s)+g(s)\Fi_\kappa(s) \ge \bigl[(1-\kappa^2)s^2+k-\mu+C_\omega\bigr]\Fi_\kappa(s),\quad \forall s\in\R.$$
Note that, by Lemma~\cite{Lemma2}, we have $1+C_\omega -\mu\le 0$, so that we can choose $\kappa_0\in(0,1)$ and
$R_0>0$ such that $(1-\kappa_0^2)R_0^2+\kappa_0+C_\omega-\mu>0$ to assure that
$$ -\Fi_{\kappa_0}''(s)+g(s)\Fi_{\kappa_0}(s)\ge 0, \quad \hbox{for}\quad \mod{s}\ge R_0.\reqlabel{PrinMax2}$$
From \cite{PrinMax1}, \cite{PrinMax2} and the Maximum Principle, we conclude that  $h(s)\le \psi_{\kappa_0}(s)$ for all
$\mod{s}\ge R_0$. So, we can choose $C(\kappa_0)\ge 1$ large enough such that
$h(s)\le C(\kappa_0)\psi_{\kappa_0}(s)$ for all $s\in\R$ to conclude the proof. \quad\cqd

\goodbreak

\smallskip 
\noindent$\bullet$ $\underline{\hbox{\sl Existence of ground states in repulsive case}}$
\smallskip 

Let us firstly consider, for $\lambda>0$, the real function $G_\lambda:\R\rightarrow\R$ defined as
$$ G_\lambda(\rho):={1\over2\lambda}\left[(1+2\lambda\rho^2)^{3/2}-(1+2\lambda\rho^2)^{1/2}\right]
       = \sqrt{1+2\lambda\rho^2}\rho^2.\reqlabel{DefGLambda}$$
It is clear that $\rho\mapsto G_\lambda(\rho)$ is a positive, convex and $C^\infty$ function such that $G_\lambda(0)=0$ and that
$$ G_\lambda'(\rho)=2f_\lambda(\rho)\rho,\quad\forall \rho\in\R.$$

Let us introduce now the functional framework which allows us to prove the existence and stability of ground states
for Eq.~\cite{OurMDE}. To do so, let $\letraEU{X}$ be the space defined by
$$\letraEU{X}:=\Bigl\{\psi\in H^1(\R)\,;\, \int_\R\left(\Mod{\psi'(s)}^2+ s^2\mod{\psi(s)}^2\right)ds<+\infty\Bigr\},\reqlabel{DefEspacoX}$$
where $H^1(\R)$ is the usual Sobolev space on $\R$.  The space
$\letraEU{X}$ is a real Hilbert space if endowed with the following usual inner product
$$(\phi|\psi)_{\letraEU{X}}:=\int_\R\left(\psi'(s)\phi'(s)+s^2\psi(s)\phi(s)\right)ds, $$
having the Hermit functions as a Hilbert orthogonal  basis.
Then, the associated norm is given by
$$\norma{\psi}{\letraEU{X}}^2:=\int_\R\left(\Mod{\psi'(s)}^2+s^2\mod{\psi(s)}^2\right)ds.$$
It is well known that the embedding of $\letraEU{X}$ into $L^p(\R)$ is compact for all $p\ge 2$
(see Proposition~6 of [\cite{KavWeis}]).

We define the {\sl energy} $E_\lambda:\letraEU{X}\rightarrow\R$ and the
{\sl charge} $Q:\letraEU{X}\rightarrow\R$ respectively by
$$\left\{\eqalign{
  E_\lambda(\psi) & := \int_\R\mod{\psi'(s)}^2ds +\int_\R s^2\mod{\psi(s)}^2ds
     +C_\omega\int_\R\sqrt{1+2\lambda\mod{\psi(s)}^2}\mod{\psi(s)}^2ds \cr
 Q(\psi) & := \int_\R\mod{\psi(s)}^2ds,\cr }\right.\reqlabel{DefEQ}$$
and denote $\Sigma_1:=\Bigl\{\psi\in\letraEU{X}\,;\, Q(\psi)=1\Bigr\}$.
By the embedding $\letraEU{X}\subset L^p(\R)$, it follows that $E_\lambda(\psi)$ is well defined because
$$ \int_\R\sqrt{1+2\lambda\psi(s)^2}\psi(s)^2\,ds\le\int_\R\left[\psi(s)^2+2\lambda\psi(s)^4\right]ds
    =\norma{\psi}{2}^2+2\lambda\norma{\psi}{4}^4.$$

With these ingredients we look for solutions  $\Fi_\mn$ of Eq.~\cite{OurMDE} that minimizes the energy
$E_\lambda$ among all functions in $\Sigma_1$.
More precisely, we look for
$$\Fi_\mn\in\Sigma_1\quad\hbox{\rm such that}\quad
     E_\lambda(\Fi_\mn)=\min\bigl\{E_\lambda(\psi)\,;\, \psi\in\Sigma_1\bigr\}.\reqlabel{VariatProb1}$$

\noindent
{\bf Theorem\ \lemlabel{ExistGrSt}:} {\sl Let $\lambda>0$ be given. Then, there exists $\Fi_\mn\in\Sigma_1$ such that }
$$E_\lambda(\Fi_\mn)=\min\Bigl\{E_\lambda(\psi)\,;\, \psi\in\Sigma_1\Bigr\}.\reqlabel{Variational}$$
\noindent
{\bf Proof:} We denote $E_\mn:=\inf\bigl\{E_\lambda(\psi)\,;\, \psi\in\Sigma_1\bigr\}$. Since
$E_\lambda(\psi)\ge0$ for all $\psi\in\letraEU{X}$,
it follows that $E_\mn\ge0$ and there exists a sequence  of minimizing functions $\{\psi_n\}_{n\in\N}$ in $\Sigma_1$, {\it i.e.},
a sequence in $\Sigma_1$ such that $\lim_{n\to+\infty} E_\lambda(\psi_n)=E_\mn$.

The fact that $G_\lambda(\rho)\ge0$ for all $\rho\in\R$ implies that
$\norma{\psi_n}{\letraEU{X}}^2\le E_\lambda(\psi_n)$ and consequently that $\{\psi_n\}_{n\in\N}$ is a bounded sequence
in $\letraEU{X}$. Therefore, it follows from the Banach-Alaoglu Theorem that there exists a subsequence
$\{\psi_{n_k}\}_{k\in\N}$ of $\{\psi_n\}_{n\in\N}$ that converges to some $\Fi_\mn$ in the weak topology
of $\letraEU{X}$, {\it i.e.}, $\psi_{n_k}\rightharpoonup\Fi_\mn$. To simplify the notation we still write $\psi_n$ for the elements of
this subsequence.

On the other hand, since $\letraEU{X}$ is compactly embedded in $L^2(\R)$, it follows that
$\Fi_\mn\in\Sigma_1$ because $\psi_n\in\Sigma_1$ for all $n$ and
$$ \lim_{n\to\infty}\int_\R\mod{\psi_n(s)}^2\,ds = \int_\R\mod{\Fi_\mn(s)}^2\,ds.$$

Since the functional $E_\lambda$ is continuous and convex, the Banach Theorem assures that $E$ is semi-continuous for the weak
topology of $\letraEU{X}$, which means that
$$ E_\lambda(\Fi_\mn)\le \liminf_{n\to+\infty} E_\lambda(\psi_n)=E_\mn$$
and consequently that $\Fi_\mn$ is a solution  \cite{VariatProb1}. \quad\cqd

\noindent
{\bf Corolary\ \lemlabel{Corol2-1}:} {\sl Each solution of \cite{VariatProb1} satisfies the equation \cite{OurMDE}.}

\noindent
{\bf Proof:}  This follows directly from the fact that $E_\lambda(\psi)$ and $Q(\psi)$ are continuously differentiable functionals
in $\letraEU{X}$. Indeed, from the Lagrange Theorem, there exists $\mu\in\R$ (a Lagrange multiplier) such that
$$E_\lambda'(\Fi_\mn)=\mu Q'(\Fi_\mn),$$
where $E_\lambda'(\psi)$ and  $Q'(\psi)$ are the Fr\'echet derivatives of $E_\lambda$ and $Q$ at $\psi\in\letraEU{X}$, respectively.
Note that this last equation is the same as \cite{OurMDE}. This completes the proof. \quad\cqd

\goodbreak

\smallskip 
\noindent$\bullet$ $\underline{\hbox{\sl Uniqueness of positive ground state}}$
\smallskip 

Let us denote by ${\cal G}_\lambda$ the set of ground states of \cite{OurMDE}, {\it i.e.},
${\cal G}_\lambda:=\bigl\{\psi\in\Sigma_1\,;\, E_\lambda(\psi)=E_\mn\bigr\}.$
It follows easily from $\norma{\psi}{\letraEU{X}}\le E_\lambda(\psi)$
that ${\cal G}_\lambda$ is a bounded set of $\letraEU{X}$ and
we have the following properties:

\smallskip\noindent
{\bf Theorem\ \lemlabel{UniqGrSt}:} {\sl Assuming $\lambda>0$, there exists a unique positive
symmetric function $\Fi_\mn\in\Sigma_1$ which is decreasing in the interval $[0,+\infty)$ and
$\Fi_\mn(0)=\max\{\Fi_\mn(s)\,;\, s\in\R\}$. Moreover}
$${\cal G}_\lambda=\bigl\{\ee{}\!^{\unim\theta}\Fi_\mn\,;\, \theta\in\R\bigr\}.$$
\goodbreak\noindent
{\bf Proof:} We proceed in two steps.

\noindent{\sl Step 1}: Let us assume for the moment that there exist two positive functions
$\Fi_1, \Fi_2\in{\cal G}_\lambda$. Then, for $0<\nu<1$, we define
$$\psi_\nu(s):=\bigl(\nu\Fi_1(s)^2+(1-\nu)\Fi_2(s)^2\bigr)^{1/2}.$$
It is clear that $\psi_\nu\in\Sigma_1$ and that
$$ \int_R s^2\mod{\psi_\nu(s)}^2ds = \nu\int_R s^2\mod{\Fi_1(s)}^2ds
    +(1-\nu)\int_R s^2\mod{\Fi_2(s)}^2ds. \reqlabel{Form1}$$
Since $\rho\mapsto G_\lambda(\rho)$ is strictly convex, it follows that
$$G_\lambda(\psi_\nu)<\nu G_\lambda(\Fi_1)+(1-\nu)G_\lambda(\Fi_2).\reqlabel{Form2}$$
Moreover, by differentiating both sides of $\psi_\nu^2=\nu\Fi_1^2+(1-\nu)\Fi_2^2$ and using
the Cauchy-Schwarz inequality, we get for each $s\in\R$,
$$\psi_\nu(s)\psi_\nu'(s)  =   \nu\Fi_1(s)\Fi_1'(s) +
    (1-\nu)\Fi_2(s)\Fi_2'(s) \le \psi_\nu(s)\sqrt{\nu\Fi_1'(s)^2+(1-\nu)\Fi_2'(s)^2},$$
from which it follows that
$$\mod{\psi_\nu'(s)}^2\le \nu\mod{\Fi_1'(s)}^2+(1-\nu)\mod{\Fi_2'(s)}^2.\reqlabel{Form3}$$
Therefore, we have from \cite{Form1}--\cite{Form3} that $E_\lambda(\psi_\nu)<\nu E_\mn+(1-\nu)E_\mn=E_\mn$,
which is impossible by the definition of $E_\mn$ and we conclude that $\Fi_1(s)=\Fi_2(s)$ for all $s\in\R$.

\smallskip
\noindent{\sl Step 2}:
We know by Theorem~\cite{ExistGrSt} that there exists at least a real function $\Fi_\mn\in{\cal G}_\lambda$.
Since $\mod{\Fi_\mn}\in\Sigma_1$ and the fact that $\letraEU{X}\subset H^1(\R)$, we have
${d\hfil\over ds}\mod{\Fi_\mn}={d\hfil\over ds}\Fi_\mn$  almost everywhere in $\R$. This implies
that $E_\lambda(\mod{\Fi_\mn})=E_\mn$ and so $\mod{\Fi_\mn}\in{\cal G}_\lambda$ is the unique positive solution of
\cite{OurMDE}.

In order to prove that $\mod{\Fi_\mn}$ is decreasing and symmetric, let $\Fi_*(s)$ be the
{\sl symmetric-decreasing rearrangement\/} of $\mod{\Fi_\mn(s)}$. As it is well known
(see [\cite{LiebLoss},\cite{Mossino}]), $\Fi_*$ is a positive and symmetric function on $\R$  such that,
for every increasing function $H:[0,+\infty)\rightarrow\R$ and every $p\ge 1$,
$$
\int_\R H\bigl(\Fi_*(s)\bigr)ds = \int_\R H\bigl(\mod{\Fi_\mn(s)}\bigr)ds\quad\hbox{\rm and}\quad
\int_\R \mod{\psi'_*(s)}^p\,ds  \le \int_\R \mod{\Fi_\mn'(s)}^p\,ds.
\reqlabel{Formula*1}$$
Hence, from \cite{Formula*1} with $H(\rho)=\rho^2$ we get that $\Fi_*\in\Sigma_1$.
On the other hand, if $c>0$,   it follows from the Hardy-Littlewood inequality,
$$\int_\R (c-s^2)^{+}\mod{\Fi_\mn(s)}^2\,ds\le \int_\R (c-s^2)^{+}\mod{\Fi_*(s)^2}\,ds,$$
which gives
$$c\int_{-\sqrt{c}}^{\sqrt c}\bigl(\mod{\Fi_\mn(s)}^2-\mod{\Fi_*(s)}^2\bigr)\,ds
   \le \int_{-\sqrt{c}}^{\sqrt c}s^2\bigl(\mod{\Fi_\mn(s)}^2-\mod{\Fi_*(s)}^2\bigr)\,ds.\reqlabel{Formula*2}$$
Therefore, using the L'Hospital rule and \cite{DecaiExp} we get for $f(s)\colon=\mod{\Fi_\mn(s)}^2-\mod{\Fi_*(s)}^2$,
$$\lim_{c\to+\infty}c\int_{-\sqrt{c}}^{\sqrt c}f(s)\,ds=-\lim_{c\to+\infty}c^2\bigl[f(\sqrt{c})-f(-\sqrt{c})\bigr]
     = -\lim_{c\to+\infty}c^2\bigl[\mod{\Fi_\mn(\sqrt{c})}^2-\mod{\Fi_\mn(-\sqrt{c})}^2\bigr]=0,$$
which implies from \cite{Formula*2} that
$$\int_\R s^2\mod{\Fi_*(s)}^2\,ds\le \int_\R s^2\mod{\Fi_\mn(s)}^2\,ds<+\infty.\reqlabel{Formula*3}$$
So, from  \cite{Formula*1} with $H=G_\lambda$, together with
the second inequality in \cite{Formula*1} with $p=2$ and \cite{Formula*3}, we get that $E_\lambda(\Fi_*)=E_\mn$, which means that
$\Fi_*\in{\cal G}_\lambda$. Therefore, from uniqueness, $\mod{\Fi_\mn(s)}=\Fi_*(s)$ for every $s\in\R$ and the proof
is complete.\quad\cqd

We will denote by $\Fi_\lambda\in {\cal G}_\lambda$ the unique real ground state
for which the minimal energy is explicitly given by
$$E_\mn=E_\lambda(\Fi_\lambda)=\int_\R\left(\Fi_\lambda'(s)^2+s^2\Fi_\lambda(s)^2\right)\!ds
    +C_\omega\int_\R \sqrt{1+2\lambda\Fi_\lambda(s)^2}\Fi_\lambda(s)^2ds.\reqlabel{EnerFiLambda}$$
The next lemma gives a characterization of the chemical potential $\mu$ as function of $\Fi_\lambda$.

\noindent
{\bf Lemma\ \lemlabel{CharcChimicPot}:} {\sl Let $\Fi_\lambda$ be the unique real function of ${\cal G}_\lambda$. Then}
$$ \mu(\Fi_\lambda)=E_\lambda(\Fi_\lambda)
    +C_\omega\lambda\int_\R{\Fi_\lambda(s)^4 \over\sqrt{1+2\lambda\Fi_\lambda(s)^2}}\,ds.\reqlabel{MuLambda}$$
\noindent
{\bf Proof:} We know that $\Fi_\lambda\in\Sigma_1$ satisfies the equation
$$-\Fi_\lambda''+s^2\Fi_\lambda+C_\omega f_\lambda\bigl(\Fi_\lambda\bigr)\Fi_\lambda=\mu(\Fi_\lambda)\Fi_\lambda.$$
By multiplying the above equation by $\Fi_\lambda$ and integrating on $\R$, we get
$$\mu(\Fi_\lambda)=\int_R\left(\Fi_\lambda'(s)^2+s^2\Fi_\lambda(s)^2\right)\!ds
  + C_\omega\int_\R f_\lambda\bigl(\Fi_\lambda(s)\bigr)\Fi_\lambda(s)^2ds.\reqlabel{DefMuFi}$$
However, it follows from \cite{DefFuncf} that
$$ f_\lambda(\rho)\rho^2  = \sqrt{1+2\lambda\rho^2}\rho^2
        +{1\over2}\left[\sqrt{1+2\lambda\rho^2}\rho^2-{\rho^2\over \sqrt{1+2\lambda\rho^2}}\right]
   = \sqrt{1+2\lambda\rho^2}\rho^2 +\lambda{\rho^4\over\sqrt{1+2\lambda\rho^2} }
       \reqlabel{GoodForm}$$
and we get \cite{MuLambda} directly from \cite{DefMuFi} and \cite{GoodForm}. \quad\cqd

\smallskip  
\noindent$\bullet$ $\underline{\hbox{\sl Stability of ground states}}$
\smallskip 

First of all, let us consider the complex extensions of $\letraEU{X}$, {\it i.e.},
$\widehat{\letraEU{X}}:=\letraEU{X}+\unim\letraEU{X}$.
Note that it can be identified with ${\letraEU{X}}\times{\letraEU{X}}$, which is a real Hilbert space if embedding with the
cartesian inner product.

In order to prove the stability of ground states of \cite{OurMDE}, we consider the Cauchy problem
$$\unim{\partial v(\tau)\over \partial\tau}={1\over 2}E_\lambda'\bigl(v(\tau)\bigr),\quad v(0,s)=v_0(s),\reqlabel{ProbCauchy}$$
where $E_\lambda'\bigl(v(\tau)\bigr)$ is the Fr\'echet derivative of $E_\lambda$ on $v(\tau)$ in $\widehat\letraEU{X}$
for all $\tau\in\R$, {\it i.e.},
$$E_\lambda'(v(\tau))=2\left(-{\partial^2 v(\tau)\over \partial s^2}+s^2v(\tau)
    +C_\omega f_\lambda\bigl(\mod{v(\tau)}\bigr)v(\tau)\right).$$

It is well known [\cite{Carles}, \cite{Oh}] that \cite{ProbCauchy} has a unique solution that is
global in time, in other words, for any $v_0\in\widehat\letraEU{X}$, there exists a unique
$v\in C\bigl([0,+\infty),\widehat\letraEU{X}\bigr)$ satisfying \cite{ProbCauchy}. In particular, if $\psi\in{\cal G}_\lambda$,
 the unique solution $u$ of \cite{ProbCauchy} such that $u(0,s)=\psi(s)$ is  the standing wave given by
$$u(\tau,s)=\ee{}\!^{-\unim\mu\tau}\psi(s).$$

Here we consider the orbital stability of ${\cal G}_\lambda$, which means that if the initial datum $v_0$ of Eq.~\cite{ProbCauchy}
is close enough to $\psi\in{\cal G}_\lambda$ in $\widehat\letraEU{X}$, then the trajectory  $v(\tau)\in\widehat\letraEU{X}$ remains
close to the set ${\cal G}_\lambda$,  as $\tau$ varies in $\R$. More precisely,

\smallskip\noindent
\noindent{\bf Definition\ \lemlabel{DefStab}:} {\sl We will say that ${\cal G}_\lambda$ is orbitally stable if, for each
$\varepsilon>0$, there exists $\delta>0$ such that if $v_0\in\widehat\letraEU{X}$ satisfies
$\displaystyle\inf_{\theta\in[0,2\pi]}\norma{v_0-\ee{}\!^{\unim\theta}\Fi_\lambda}{\widehat\letraEU{X}}<\delta$,
the solution of \cite{ProbCauchy} satisfies}
$$\sup_{\tau\in\R}\left(\inf_{\theta\in[0,2\pi]}\norma{v(\tau)
     -\ee{}\!^{-\unim\mu(\Fi_\lambda)\tau}\ee{}\!^{\unim\theta}\psi}{\widehat\letraEU{X}}\right)<\varepsilon.$$

In order to prove the stability of ${\cal G}_\lambda$, we follow Cazenave-Lions [\cite{CazeLions}], where
the orbital stability is a direct consequence of the well known conservation laws
that hold for all solutions of \cite{ProbCauchy}: {\sl If $v(\tau,s)$ is a solution of \cite{ProbCauchy}
such that $v(0,s)=v_0(s)$, then, for each $\tau\in\R$, we have}
$$Q\bigl(v(\tau)\bigr) = Q(v_0)\quad\hbox{\rm and}\quad E_\lambda\bigl(v(\tau)\bigr) = E_\lambda(v_0).$$
The above identities are known respectively as {\sl conservation of charge\/} and {\sl conservation of energy\/} and are
easily obtained by multiplying the equation respectively by the complex conjugates of $-\unim v(\tau)$ and
$\partial v(\tau)/\partial\tau$.

\smallskip\goodbreak
\noindent{\bf Theorem\ \lemlabel{Thm2}:} {\sl Assuming $\lambda>0$, the set ${\cal G}_\lambda$ is stable
in the sense of the Definition~\cite{DefStab}.}

\smallskip\noindent
{\bf Proof:}  We argue by contradiction.
If ${\cal G}_\lambda$ is not stable in the sense of Definition \cite{DefStab}, there exists
$\varepsilon_0>0$ such that for all $n\in\N$, we can find  $v_{0n}\in\widehat\letraEU{X}$ satisfying
$$r_n:=\inf_{\theta\in[0,2\pi]}\norma{v_{0n}-
    \ee{}\!^{\unim\theta}\Fi_\lambda}{\widehat\letraEU{X}}<{1\over n}\reqlabel{PrimaCond}$$
and
$$\sup_{\tau\in\R}\left(\inf_{\theta\in[0,2\pi]}\norma{v_n(\tau)
      -\ee{}\!^{-\unim\mu_\lambda\tau}\ee{}\!^{\unim\theta}\Fi_\lambda}{\widehat\letraEU{X}}\right)\ge
      \varepsilon_0,\reqlabel{DuaCond}$$
where $v_n\in C(\R;\widehat\letraEU{X})$ is the unique solution of \cite{ProbCauchy} with initial datum $v_{0n}$
and $\mu_\lambda:=\mu(\Fi_\lambda)$.

Let $\theta_n\in[0,2\pi]$ such that
$ r_n\le \norma{v_{0n}-\ee{}\!^{\unim\theta_n}\Fi_\lambda}{\widehat\letraEU{X}}<1/n$.
By compactness, there exist $\theta_*\in[0,2\pi]$ and a subsequence of $\{\theta_n\}_{n\in\N}$ (by simplicity still denoted
by  $\{\theta_n\}_{n\in\N}$),  such that   $\theta_n\rightarrow\theta_*$, from which we may infer that
$$
   \lim_{n\to\infty}\norma{ v_{0n}-\ee{}\!^{\unim\theta_*}\Fi_\lambda}{\widehat\letraEU{X}}=0.
$$
Moreover, we can select  a sequence $\{\tau_n\}_{n\in\N}$
in $\R$ such that
$$ \inf_{\theta\in[0,2\pi]} \norma{v_n(\tau_n)
   -\ee{}\!^{-\unim\mu_\lambda\tau_n}\ee{}\!^{\unim\theta}\Fi_\lambda}{\widehat\letraEU{X}}\ge\varepsilon_0/2,
   \quad\forall n\in\N.   \reqlabel{QuadCond}$$
Now, denoting $\psi_n:=\ee{}\!^{\unim\mu_\lambda\tau_n}v_n(\tau_n)$, we have from \cite{QuadCond}
$$\inf_{\theta\in[0,2\pi]}\norma{\psi_n-\ee{}\!^{\unim\theta}\Fi_\lambda}{\widehat\letraEU{X}}\ge\varepsilon_0/2,
\quad\forall n\in\N.    \reqlabel{QuintCond}$$

From the conservation laws, we get, for $n\to+\infty$,
$$ Q(\psi_n) = Q(v_{0n}) \rightarrow   Q(\ee\!^{\unim\theta_*}\Fi_\lambda)=1 \quad\hbox{and}\quad
      E_\lambda(\psi_n) = E_\lambda(v_{0n}) \rightarrow   E_\lambda(\ee\!^{\unim\theta_*}\Fi_\lambda)=E_\mn.$$

Note that the sequence $\{E_\lambda(\psi_n)\}_{n\in\N}$ is convergent and as
$\norma{\psi_n}{\widehat\letraEU{X}}\le E_\lambda(\psi_n)$,
there exists a subsequence still indexed by $n$ that converges to some $\psi_*$ in the weak
topology of $\widehat\letraEU{X}$, {\it i.e.}, $\psi_n\rightharpoonup\psi_*$ in $\widehat\letraEU{X}$.
By remembering that $\widehat\letraEU{X}$ is compactly embedded in $L^p(\R)$ for all $p\ge2$,
it follows that  $\psi_n\rightarrow\psi_*$ in $L^2(\R)$ and we have $\psi_*\in\Sigma_1$.

On the other hand, the functional $J(\psi):=\int_\R G_\lambda\bigl(\psi(s)\bigr)\,ds$ is convex and well defined in
$V:=L^2(\R)\cap L^4(\R)$, which is a Banach space  for the norm $\norma{\ }{V}:=\norma{\ }{2}+\norma{\ }{4}$.
So, $J$ is continuous in $V$ and then $J(\psi_n)\rightarrow J(\psi_*)$ as $n\to+\infty$.

Since
$$ \norma{\psi_*}{\widehat\letraEU{X}}^2\le \liminf_{n\to\infty}\norma{\psi_n}{\widehat\letraEU{X}}^2
  = \lim_{n\to\infty}\bigl(E_\lambda(\psi_n)-C_\omega J(\psi_n)\bigr)
  = E_\mn-C_\omega J(\psi_*)$$
and  $\psi_*\in\Sigma_1$, it follows that $E_\lambda(\psi_*)=E_\mn$, which means that  $\psi_*\in\cal{G}_\lambda$
and so, $\displaystyle \inf_{\theta\in[0,2\pi]}\norma{\psi_*-\ee\!^{\unim\theta}\Fi_\lambda}{\widehat\letraEU{X}} =0$.
This is in contratiction of \cite{QuintCond} because
$$ \lim_{n\to\infty}\norma{\psi_n}{\widehat\letraEU{X}}^2
   = \lim_{n\to\infty}\bigl(E_\lambda(\psi_n)-C_\omega J(\psi_n)\bigr)
   = E_\mn-C_\omega J(\psi_*)=\norma{\psi_*}{\widehat\letraEU{X}}^2.$$
which implies that $\psi_n\rightarrow \psi_*$ in $\widehat\letraEU{X}$ and consequently
$$ \norma{\psi_*-\ee\!^{\unim\theta}\Fi_\lambda}{\widehat\letraEU{X}}
   =\lim_{n\to\infty}\norma{\psi_n-\ee\!^{\unim\theta}\Fi_\lambda}{\widehat\letraEU{X}}
      \ge \eps_0/2,\quad\forall\theta\in[0.2\pi].\quad\cqd$$
%
\noindent{\bf\numsection.\ A  relation between chemical potential and energy as functions of $\lambda>0$}\par
\smallskip 


In this section we report the relations between the energy $E_\lambda(\Fi_\lambda)$ and the chemical potential
$\mu(\Fi_\lambda)$ as functions of $\lambda$, which are well defined from the uniqueness of
$\Fi_\lambda$, for all $\lambda\ge0$. More precisely, let $E_\mn,\mu_\mn:[0,+\infty)\rightarrow\R$ be the functions
defined by
$$E_\mn(\lambda):= E_\lambda(\Fi_\lambda)\quad\hbox{\rm and}\quad \mu_\mn(\lambda):=\mu(\Fi_\lambda).$$

\goodbreak
\noindent{\bf Proposition\ \lemlabel{PropELambda}:} {\sl The function $E_\mn$ is a positive, increasing
and concave function such that
$$\lim_{\lambda\to+\infty}E_\mn(\lambda)=+\infty.\reqlabel{LimEtoInfty}$$
Assuming that  the curve
$\lambda\mapsto\Fi_\lambda\in\letraEU{X}$ is  continuously differentiable in $(0,+\infty)$, we have
$$ E_\mn \in C^1(0,+\infty)\quad\hbox{and}\quad  {d\hfil\over d\lambda}E_\mn(\lambda)
    = C_\omega\int_\R {\Fi_\lambda(s)^4\over \sqrt{1+2\lambda\Fi_\lambda(s)^2}}\,ds.\reqlabel{dEdlambda}$$
Moreover, $\mu_\mn$ is continuous, satisfy the same limit \cite{LimEtoInfty}
and is related to $E_\mn$ by }
$$E_\mn(\lambda)={1\over\lambda}\int_0^\lambda \mu_\mn(\xi)\,d\xi.\reqlabel{MyFormula} $$
\relax \noindent
{\bf Proof:} For each $\psi\in\Sigma_1$ and $s\in\R$, the map $\lambda\mapsto G_\lambda\bigl(\psi(s)\bigr)$ is
increasing and concave, which imply that $\lambda\mapsto E_\lambda(\psi)$  is also increasing and concave.
As it is well known that the minimum of a family of increasing and concave functions is itself increasing and concave,
the same is true for $E_\mn(\lambda)$.

To prove the limit \cite{LimEtoInfty}, let us proceed by contradiction. Since $E_\mn$ is increasing, if \cite{LimEtoInfty}
does not hold, we have
$\lim_{\lambda\to+\infty}E_\mn(\lambda)=L$ for some $L>0$. In particular, for $\lambda_n:=n$,
the sequence $\{\Fi_n\}_{n\in\N}$ is bounded in \letraEU{X}. Hence, there exists $\Fi_\infty\in\letraEU{X}$ such that
$\Fi_n\rightharpoonup \Fi_\infty$ in $\letraEU{X}$, and by the compact embedding, $\Fi_\infty\in L^p(\R)$ for all $p\ge2$.
So, $\Fi_\infty\in\Sigma_1$ and
$$ L\ge E_\mn(n) \ge \int_\R\sqrt{1+2n\Fi_n(s)^2}\Fi_n(s)^2ds\ge \sqrt{2n}\norma{\Fi_n}{3}^3\quad  \forall n\in\N,$$
which is impossible.

From the definition of the minimal energy, we have
$${d\hfil\over d\lambda}E_\mn(\lambda) = 2\int_\R\Fi_\lambda'(s) {d\hfil\over d\lambda}\Fi_\lambda'(s)\,ds
    + 2\int_\R s^2\Fi_\lambda(s) {d\hfil\over d\lambda}\Fi_\lambda(s)\,ds + C_\omega {d\hfil\over d\lambda}H(\lambda),
    \reqlabel{My1}$$
where $ H(\lambda):= \int_\R\sqrt{1+2\lambda\Fi_\lambda(s)^2}\Fi_\lambda(s)^2ds $.
By  direct calculation we get
$$ {d\hfil\over d\lambda}H(\lambda) =
   2\int_\R{1+3\lambda\Fi_\lambda(s)^2 \over \sqrt{1+2\lambda\Fi_\lambda(s)^2}}\Fi_\lambda(s)
     {d\hfil\over d\lambda}\Fi_\lambda(s)\,ds
      + \int_{\R}{\Fi_\lambda(s)^4 \over \sqrt{1+2\lambda\Fi_\lambda(s)^2}}\,ds.\reqlabel{My2} $$
Now, remembering that
$ E'(\psi) = -2\psi''(s)+2s^2\psi(s)+2C_\omega\displaystyle{1+3\lambda\psi(s)^2 \over \sqrt{1+2\lambda\psi(s)^2}}\psi(s)$,
it follows from \cite{My1} and \cite{My2}  that
$$ {d\hfil\over d\lambda}E_\mn(\lambda)
   =\left<E'(\Fi_\lambda):{d\hfil\over d\lambda}\Fi_\lambda\right>
       + C_\omega\int_\R{\Fi_\lambda(s)^4\over \sqrt{1+2\lambda\Fi_\lambda(s)^2}}\,ds,$$
where $\Bigl<\, :\,\Bigr>$ denotes the duality product in  $\letraEU{X}'\times\letraEU{X}$.
Since $E'(\Fi_\lambda)=\mu(\lambda)Q'(\Fi_\lambda)$ and $Q(\Fi_\lambda)=1$ for all $\lambda\ge0$,
we have
$$ \left<E'(\Fi_\lambda):{d\hfil\over d\lambda}\Fi_\lambda\right>
       =\mu(\lambda)\left<Q'(\Fi_\lambda):{d\hfil\over d\lambda}\Fi_\lambda\right>
       = \mu(\lambda){d\hfil\over d\lambda}Q(\Fi_\lambda)=0,$$
from which \cite{dEdlambda} is proved.

On the other hand, by Lemma~\cite{CharcChimicPot}, it follows that $\mu_\mn(\lambda)\ge E_\mn(\lambda)$ for all $\lambda>0$,
which implies that $\mu_\mn$ has the same limit as $E_\mn$ at infinity.  Furthermore, from \cite{MuLambda}
and \cite{dEdlambda} we get
$$\mu_\mn(\lambda)=E_\mn(\lambda)+\lambda {d\hfil\over d\lambda}E_\mn(\lambda)
   ={d\hfil\over d\lambda}\bigl(\lambda E_\lambda(\lambda)\bigr),$$
from which we obtain \cite{MyFormula} after integrating the above identity on $[0,\lambda]$.
\quad\cqd

\smallskip\goodbreak
\noindent{\bf Remark\ \lemlabel{ObsLambNeg}:} The parameter $\lambda$ in our context is related
to the $s$-wave scattering length $a$ (or the product $aN$) and, as inferred in [\cite{MunozDelgado}], the model \cite{MDE}
is also valid for the attractive case ($a<0$).  So, it would be important to assure this fact from a
mathematical point of view, which means to prove that \cite{OurMDE} admit solutions for $\lambda<0$.
Therefore, from the definition of the energy, one must be  assured that
$2\lambda\Fi_\lambda(s)^2+1\ge0$ for all $s\in\R$ and $\lambda$ in a certain interval $(-\lambda_*,0)$
for some $\lambda_*>0$,  or
from Theorem~\cite{UniqGrSt}, $2\lambda\Fi_\lambda(0)^2+1\ge0$ for $\lambda>-\lambda_*$.

\medskip 
\noindent{\bf\numsection.\ Approximation formulae}\par
\smallskip 

Based on the exponential decay as proved in Theorem~\cite{ExpDecay} and the properties
presented in Theorem~\cite{UniqGrSt}, we are led to consider the curve $\Gamma_1\subset\Sigma_1$ defined by
$$\Gamma_1:=\bigl\{\psi_\kappa(s):=(\kappa/\pi)^{1/4}\exp(-\kappa s^2/2), \quad \kappa>0\bigr\}\reqlabel{DefFiKappa}$$
as trial functions to evaluate the minimal energy $E_\mn(\lambda)$.

Note that if $\lambda=0$, it follows from \cite{DefEQ} that
$E_0(\psi)=\norma{\psi}{\letraEU{X}}^2+C_\omega\norma{\psi}{2}^2$, which reaches its minimum
at $\Fi_1(s)=(1/\pi)^{1/4}\ee^{-s^2/2}$ and $E_\lambda(\Fi_1)=\mu(\Fi_1)=1+C_\omega$.

A direct calculation of the energy at $\psi_\kappa\in\Gamma_1$ gives
$$E_\lambda(\psi_\kappa) = {1\over 2}\left(\kappa+{1\over\kappa}\right)
        +C_\omega\int_\R\sqrt{1+2\lambda\psi_\kappa(s)^2}\psi_\kappa(s)^2ds.\reqlabel{EappGeral}$$

By the same argument in the proof of Theorem~\cite{ExistGrSt} we can show that there exists a
unique $\kappa(\lambda)>0$ such that
$$ E_\app(\lambda):=E_\lambda(\psi_{\kappa(\lambda)}) = \min\bigl\{E_\lambda(\psi_\kappa)\,;\, k>0\bigr\}.\reqlabel{ProbVariatApp}$$
Unfortunately there is no explicit solution of the integral in \cite{EappGeral} to calculate $\kappa(\lambda)$ from this variational problem
involving a simple real function of $\kappa$.

The direct calculation of the Taylor series of $G_\lambda(\rho)$ definied in \cite{DefGLambda} gives
$$G_\lambda(\rho) = {1\over 2\lambda}\left[2\lambda\rho^2
   +\sum_{n=2}^\infty {(-1)^n\over n!}{(2n-4)!\over(n-2)!}{2n \over 2^{2n-2}}\bigl(2\lambda\rho^2\bigr)^n\right],$$
which converges absolutely and uniformly for $2\lambda\rho^2$ in compact subsets of $(-1,1)$.

For $\rho=\psi_\kappa(s)$, we have
$$ E_\lambda(\psi_\kappa)  =  {1\over 2}\left(\kappa+{1\over\kappa}\right) + C_\omega\left(1
     +\sum_{n=2}^\infty a_n\lambda^{n-1}\kappa^{(n-1)/2}\right),\quad
        a_n=\displaystyle{(-1)^n \over n! }{(2n-4)!\over(n-2)!}{\sqrt{n}\over2^{n-2}}
        \left({1\over\pi}\right)^{(n-1)/2}, \reqlabel{TaylorKappa}$$
which is convergent as long as $\kappa<\pi/4\lambda^2$.

This development up to order $n=2$ allows us to recover the energy of a cubic model already widely
discussed in [\cite{CipGondTrall}]. Moreover, if we consider the Taylor's development of $G_\lambda(\rho)$
up to order $n=3$, we get the energy of a cubic-quintic  model (see  [\cite{TrallCipLiew}]).

To illustrate our numerical results, we assume that $C_\omega=1$. Then, from \cite{TaylorKappa}, we have
$$ {d\hfil\over d\kappa}E_\lambda(\psi_\kappa) = {1\over 2}\left(1-{1\over\kappa^2}\right) +
   \sum_{n=2}^\infty {n-1\over2}a_n\lambda^{n-1}\kappa^{(n-3)/2}=0,\reqlabel{Eq.2}$$
or equivalently,
$$ 2\kappa^2 {d\hfil\over d\kappa}E_\lambda(\psi_\kappa)=\kappa_2-1 +
\sum_{n=2}^\infty  (n-1)a_n\lambda^{n-1}\kappa^{(n+1)/2}=0.$$
By considering Taylor's approximation up to orden $n=6$, we get
$$ 2\kappa^2{d\hfil\over d\kappa}E_\lambda(\psi_\kappa) = 0 \iff f(\lambda,\kappa)=1,\reqlabel{Eq.FLKappa}$$
for
$ f(\lambda,\kappa):=a_2\lambda\kappa^{3/2} + (1+2a_3\lambda^2)\kappa^2 + 3a_4\lambda^3\kappa^{5/2}
       + 4a_5\lambda^4\kappa^3 + 5a_6\lambda^3\kappa^{7/2}$, where
$$ a_2={1\over\sqrt{2\pi}},\quad
     a_3=-{1\over2\sqrt{3}\pi},\quad
     a_4={1\over4\pi\sqrt\pi},\quad
     a_5=-{\sqrt5\over8\pi^2},\quad
     a_6={7\sqrt6\over48\pi^2\sqrt\pi}.\reqlabel{Coefic}$$

In the sequel, we denote  $E_\app(\lambda)$ and $\mu_\app(\lambda)$ respectively the energy and the
corresponding chemical potential, with approximation up to order $n=6$. We present their graphics based on
determining the solution $\kappa(\lambda)$ of Eq.~\cite{Eq.FLKappa}, which are unique for any $\lambda\ge0$
(even for $\lambda<0$ if $\mod\lambda$ is not very large). The first figure shows the graphic of the
function $\lambda\mapsto\kappa(\lambda)$.

\setbox11=\hbox{\epsfxsize=7cm\epsfbox{graficos_BEC_n6.1}}
$$\hbox{\box11}$$
{\parindent=1.5truecm
\narrower\noindent
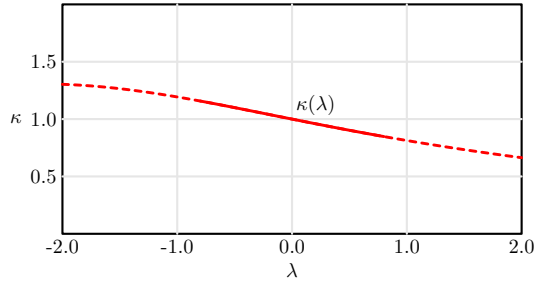
{\bf Fig.~1:}  {\letranota Graphic of $\kappa(\lambda)$ for $\lambda\in[-2,2]$.
The solid line corresponds to Eq.~\cite{Eq.FLKappa} in the interval of convergence of \cite{TaylorKappa}. }
\par }

\smallskip
By using the values of $\kappa(\lambda)$ in the Eq.~\cite{EappGeral} and numerical integration of the non-quadratic term
with $\psi_{\kappa(\lambda)}$,  we obtain the graphics of $E_\app(\lambda)$ for $\lambda\in(-\lambda_\app,2]$,
where $\lambda_\app\approx 0.78$. Similarly, by using the values of $E_\app(\lambda)$ and numerical integration
in Eq.~\cite{MuLambda}, we obtain the graphic of $\mu_\app(\lambda)$.

\setbox11=\hbox{\epsfxsize=7cm\epsfbox{graficos_BEC_n6.7}}
$$\hbox{\box11}$$
{\parindent=1.5truecm
\narrower\noindent
{\bf Fig.~2:} {\letranota Graphics of  $E_\app(\lambda)$  (blue) and $\mu_\app(\lambda)$ (purple) for
$\lambda\in(-\lambda_\app,2]$. In both cases, the solid lines correspond to the interval of convergence
of \cite{TaylorKappa}}\par }

\noindent
\noindent{\bf Remark\ \lemlabel{DerivIMplicMeth}:} We could calculate the function $\kappa(\lambda)$
by considering the implicity derivation of Eq.~\cite{Eq.FLKappa}. More precisely,
$$ {d\hfil\over d\lambda}f\bigl(\lambda,\kappa(\lambda)\bigr)
    ={\partial f\over\partial\lambda}+{\partial f\over\partial\kappa}{d\kappa\over d\lambda}=0,\reqlabel{ImplDeriv}$$
which allows us to determine the function $\kappa(\lambda)$ as the solution of an initial value problem for ODE, namelly,
$$ {d\kappa\over d\lambda}=
  -{\partial f\over \partial\lambda}\left({\partial f\over \partial\kappa}\right)^{-1},\quad \kappa(0)=1.$$
This method produces similar results as before if $\lambda$ is restrict to the interval of convergence of \cite{TaylorKappa}.

\medskip 
\noindent{\bf\numsection.\ Generalized Thomas-Fermi Approximation}\par
\smallskip 

Under the condition that the kinetic energy be small enough, {\it i.e.} if
$$ \int_\R\mod{\Fi'(s)^2}\,ds\ll C_\omega\int_\R f_\lambda\bigl(\Fi(s)\bigr)\Fi(s)^2\,ds,$$
the term $\Fi''$  is negligible, which means that, for $\lambda>0$ large enough,
the equation  is reduced to
$$ f_\lambda(\Fi(s))={\mu-s^2\over C_\omega},\reqlabel{Thomas-Fermi1}$$
where $ f_\lambda(\rho)={3\over 2}\sqrt{1+2\lambda\rho^2}-{1\over 2\sqrt{1+2\lambda\rho^2}}$.

It is clear that $\rho\mapsto f_\lambda(\rho)$ is strictly increasing in $[0,+\infty)$
and its inverse can be easily calculated. Indeed, from its definition we can write
$f_\lambda(\rho) := h\bigl(\sqrt{1+2\lambda\rho^2}\bigr)$, where $\xi:=h(\eta)={1\over 2}\bigl(3\eta-{1\over \eta}\bigr)$,
whose inverse is
$$h^{-1}(\xi)={1\over 3}\left(\xi+\sqrt{\xi^2+3}\right).$$
Therefore, from \cite{Thomas-Fermi1} we get
$$ \sqrt{1+2\lambda\Fi(s)^2}  = h^{-1}\left({\mu-s^2\over C_\omega}\right)
          ={\mu\over 3C_\omega}\left(1-{s^2\over\mu}
              +\sqrt{\left(1-{s^2\over\mu}\right)^2+{3C_\omega^2\over\mu^2}}\right),$$
or equivalently,
$$ 1+2\lambda\Fi(s)^2 = {2\mu^2\over 9C_\omega^2}
 \left[\left(1-{s^2\over\mu}\right)^2
   +  \left(1-{s^2\over\mu}\right)\sqrt{ \left(1-{s^2\over\mu}\right)^2+{3C_\omega^2\over\mu^2}}\,
  \right]+{1\over3}
$$
From the fact that $\mu$ is large if $\lambda>0$ is large enough (see Proposition~\cite{PropELambda}),
we can assume that $3C_\omega^2/\mu^2$ is negligible to get the following (first) approximation
$$ 1+2\lambda\Fi(s)^2  = {4\mu^2\over 9C_\omega^2}\left(1-{s^2\over\mu}\right)^2+{1\over3}.\reqlabel{FirstAprox}$$
Hence, considering that $\int_{-\sqrt\mu}^{\sqrt\mu}\Fi_\lambda(s)^2ds\approx1$
(and remembering the symmetry of $\Fi$), we get
$$ {4\mu^2\sqrt\mu\over 9C_\omega^2}\int_0^1(1-u^2)^2du-{2\sqrt\mu\over3}=\lambda.\reqlabel{FirstMuTF}$$
On the other hand, instead of simply discard the term $3C_\omega/\mu^2$, we can consider the approximation of
$\sqrt{x+h}$ as $\sqrt{x}+h/2\sqrt{x}$, with $x=(1-s^2/\mu)^2$ and $h=3C_\omega^2/\mu^2$, to get
$$ 1+2\lambda\Fi(s)^2  = {4\mu^2\over 9C_\omega^2}\left(1-{s^2\over\mu}\right)^2+{2\over3},\reqlabel{SecondAprox}$$
which gives the (second) approximation
$$ {4\mu^2\sqrt\mu\over 9C_\omega^2}\int_0^1(1-u^2)^2du-{\sqrt\mu\over3}=\lambda.\reqlabel{SecondMuTF}$$
The solutions of \cite{FirstMuTF} and \cite{SecondMuTF} provide good approximations for $\mu$ as functions of $\lambda$,
as shown in the figure below (where we have considered $C_\omega=1$  for numerical calculations).
\setbox11=\hbox{\epsfxsize=7cm\epsfbox{Test_Thomas_Fermi.1}}
$$\hbox{\box11}$$
{\parindent=1.5truecm
\narrower\noindent
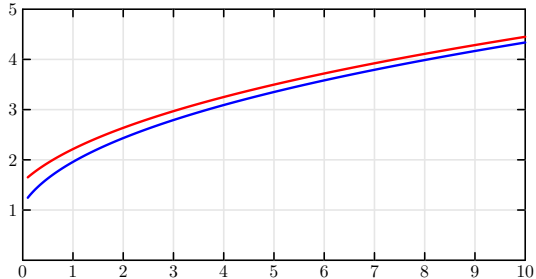
{\bf Fig.~3:} {\letranota Graphics of  $\mu(\lambda)$ for $\lambda\in(0,10]$ by Thomas-Fermi aproximation.
The red line corresponds to Eq.~\cite{FirstMuTF} and the blue one, to Eq.~\cite{SecondMuTF}.
\par} }

\goodbreak
\medskip 
\noindent{\bf\numsection.\ Conclusions}\par
\smallskip 

In this paper we present some important mathematical  properties of the ground-state solutions of
one-dimensional effective models proposed by  A.~Mu\~noz Mateo and V.~Delgado to describe the cigar-shaped
Bose-Einstein condensates in harmonic trap potentials. In Section~2, we prove, in the case of
repulsive interatomic forces ($\lambda>0$), that there exists a unique positive and symmetric ground state
$\Fi_\mn$, decreasing for $s>0$ (Theorems \cite{ExistGrSt}  and \cite{UniqGrSt}), which is orbitally stable
in the Hilbert space of finite energy functions (Theorem~\cite{Thm2}). Moreover, like any other solution of
Eq.~\cite{OurMDE}, such ground state has a Gausssian-like exponential asymptotic decay (Theorem~\cite{ExpDecay}).
In section~3 we describe the minimal energy $E_\mn$ and the corresponding chemical potential $\mu_\mn$ as functions
of the parameter $\lambda$, for which  some general properties are shown (Proposition~\cite{PropELambda}).
All theses properties suggest that Gaussian functions could be used as approximations of those BECs. Based on this idea,
we present in Section~4 some numerical examples. Finally, in Section~5, we describe an extension of the Thomas-Fermi
method to calculate the chemical potential in the context of large $\lambda>0$.

\bigskip 
\noindent{\bf\numsection.\  References}\par
\kern-.6truecm
\MakeBibliography{}

\bye